\documentclass[12pt]{article}

\usepackage{amsthm,amsmath,amssymb}

\usepackage{graphicx}
 \usepackage{epsfig}
  \usepackage{latexsym, epsfig, psfrag,eepic,colordvi,bm}
  \usepackage{graphics,color}
  \usepackage{graphics}
  \usepackage{graphicx}
  \usepackage{epsfig}
  \usepackage[all]{xy}
\usepackage{cancel}
  \usepackage{amssymb}
  \usepackage{amsthm,amsmath}

  \usepackage{latexsym, epsfig, psfrag,eepic,colordvi,bm}
  \usepackage{graphics,color}
  \usepackage{amsmath,amsfonts,amssymb,amscd}
  \usepackage[all]{xy}
  \usepackage{epstopdf}

  \usepackage{amsthm,amsmath,amssymb}

  \usepackage{graphicx}
  \usepackage{amsthm,amsmath,amssymb}

  \usepackage{graphicx,epsfig}
\usepackage[colorlinks=true,citecolor=black,linkcolor=black,urlcolor=blue]{hyperref}

\usepackage[all]{xy}
\usepackage{graphicx}
\usepackage{mathptmx}

\usepackage{latexsym, epsfig,cite, psfrag,eepic,colordvi,bm}
\usepackage{graphics,color,graphicx}
\usepackage[justification=centering]{caption}
\usepackage{amsmath,amsfonts,amssymb,amscd}
\usepackage[all]{xy}

\theoremstyle{plain}
\newtheorem{theorem}{Theorem}[section]

\newtheorem{lemma}{Lemma}[section]
\newtheorem{corollary}{Corollary}[section]
\newtheorem{proposition}{Proposition}[section]

\theoremstyle{definition}

\theoremstyle{remark}

\newcommand{\genstirlingI}[3]{%
	\genfrac{[}{]}{0pt}{#1}{#2}{#3}%
}

\newcommand{\stirlingI}[2]{\genstirlingI{}{#1}{#2}}

\date{}

\title{The Harer-Zagier and Jackson formulas and new results for one-face bipartite maps}

\author{Ricky X. F. Chen~\footnote{ORCID: 0000-0003-1061-3049}\\
	\small School of Mathematics, Hefei University of Technology \\[-0.8ex]
	\small Hefei, Anhui 230601, P.~R.~China\\[-0.8ex]
	\small\tt xiaofengchen@hfut.edu.cn
}

\begin{document}

\maketitle

\begin{abstract}

The study of bipartite maps (or Grothendieck's dessins d'enfants) is closely connected with
geometry, mathematical physics and free probability.
Here we study these objects from their permutation factorization formulation using a novel
character theory approach. We first present some general symmetric function expressions for the number
of products of two permutations respectively from two arbitrary, but fixed, conjugacy classes indexed by $\alpha$ and $\gamma$
which produce a permutation with $m$ cycles.
Our next objective is to derive explicit formulas for the cases where
$\alpha$ corresponds to full cycles, i.e., one-face bipartite maps.
We prove a far-reaching explicit formula, and
show that the number for any $\gamma$ can be iteratively reduced to that of products of two full cycles,
which implies an efficient dimension-reduction algorithm for building a database of all these numbers.
Note that the number for products of two full cycles can be computed by the Zagier-Stanley formula.
Also, in a unified way,
we easily prove the celebrated Harer-Zagier formula
and Jackson's formula, and we obtain explicit formulas for
several new families as well.


\end{abstract}
\section{Introduction}\label{sec1}

The study of maps is ubiquitous, as maps are closely connected with
geometry, mathematical physics, free probability and algebra.
 A map of genus $g$ is a $2$-cell embedding of a connected graph (loops and multi-edges allowed) in an orientable topological surface
 of genus $g$
 such that all faces are homomorphic to an open disk.
 The genus $g$ of the topological surface relates to the number $v$ of vertices,
 $e$ of edges and $f$ of faces via the Euler characteristic formula: 
 $$
 v-e+f=2-2g.
 $$
 A bipartite map is a map where there exists a way of coloring the vertices such that
two adjacent vertices have different colors and at most two colors in total are used.
Obviously, the underlying graph of a bipartite map has no loops.
 Bipartite
maps are also known as Grothendieck’s dessins d’enfants in studying Riemann surfaces~\cite{Grothendieck,zap,lan-zvon},
and equivalent to hypermaps~\cite{walsh2}.

\begin{figure}[!htb]
	\centering
	\includegraphics[scale=.8]{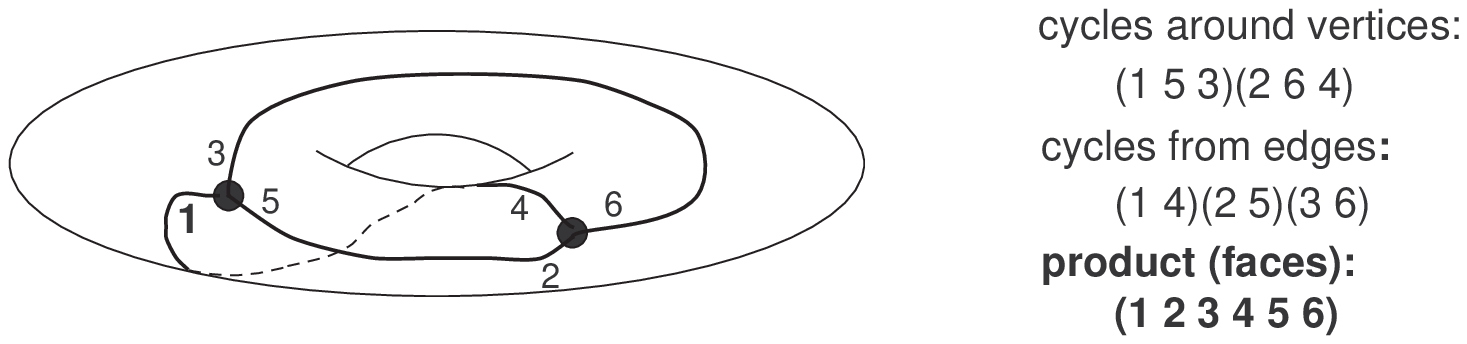}
	\caption{A map (rooted at $1$) of genus one with two vertices, three edges and one face.}\label{fig:map}
\end{figure}

An end of an edge is called a half edge or a dart, so there are $2n$ half edges in a graph with $n$
edges.
It is long known that an embedding is completely determined by the cyclic orders of the half edges
around the vertices (see e.g., Edmonds~\cite{edmonds}).
Consequently, there is a well established algebraic formulation of maps which works as follows.
If we assign the numbers in $[2n]=\{1,2,\ldots, 2n\}$ to the half edges, then two permutations on $[2n]$ immediately arise:
the cycles of one permutation $\sigma_1$ come from the half edges around (counterclockwisely) the vertices of the graph and the other, $\sigma_2$, consists of $n$ cycles
where the two half edges of an edge give a length-two cycle.
The product $\omega=\sigma_2\sigma_1$ (multiplying from right to left) of these two permutations actually encode the faces of the embedding.
As such, enumeration of maps can be translated into permutation product problems.
The genus $g$ of the underlying surface is then encoded in the numbers of cycles of the involved permutations.
See Figure~\ref{fig:map} for an example.
One hidden requirement on the $\sigma$-pairs is that the group generated by the two permutations acts transitively
on the set $[2n]$, which is implied by the connectivity of the underlying graph.

If one of the half edges is particularly distinguished, the map is called rooted.
All considered maps (and bipartite maps) are rooted in this work, and rooted at half edge $1$ if not specified otherwise. 
Two maps are viewed equivalent if there is an orientation preserving homeomorphism
of the underlying surface such that vertices are mapped to vertices, edges are mapped
to edges, and in particular, the root half edge of one map is mapped to the root half edge
of the other. Algebraically, two maps $(\sigma_1, \sigma_2, \omega)$
and $(\sigma_1', \sigma_2',\omega')$
are equivalent if there exists a permutation $\pi$ such that
\begin{align*}
	\sigma_1'=\pi \sigma_1 \pi^{-1}, \quad \sigma_2'=\pi \sigma_2 \pi^{-1}, \quad \omega'=\pi \omega \pi^{-1}, \quad \pi(1)=1.
\end{align*}

Although bipartite maps of $n$ edges can be also encoded into triples of permutations like maps described above,
conversely, the map corresponding to a triple of permutations may not be bipartite.
To remedy this, we employ a slightly different encoding.
That is, given a rooted bipartite map of $n$ edges, we first assign white and black colors
to the vertices such that the vertex incident to the root is in black, and label the edges $1,2,\ldots, n$. We then encode it into a triple of permutations on the
set $[n]$ as follows: the cycles of $\sigma_1$ come from the edges around the black vertices,
the cycles of $\sigma_2$ come from the edges around the white vertices, and $\omega= \sigma_2 \sigma_1$ gives the faces.
Note that the length of a cycle of $\sigma_2$ is no longer necessarily two.
Equivalent bipartite maps are defined analogously.
As such, enumerating bipartite maps is tantamount to the enumeration of products of permutations satisfying certain
conditions.

The enumeration of maps may date back to Tutte in the 1960s in a series of papers about planar maps (i.e., genus zero maps).
Explicit formulas for maps of arbitrary genus were first obtained by Walsh and Lehman~\cite{walsh1}.
In the study of moduli space of curves,
the celebrated Harer-Zagier formula and recurrence satisfied by the numbers of
maps with only one face were obtained via integral over matrices~\cite{harer-zagier}. Making use of the permutation product model,
Jackson~\cite{jac} first provided an approach to compute these numbers using the symmetric group
character theory, and Zagier~\cite{zag} later also provided a group character approach but shorter.   
Very recently, Pittel~\cite{pittel} gave another derivation of the Harer-Zagier formula via group characters
in the form of Fourier transform.
Combinatorial approaches attacking this topic have also
received lots of interests. We refer to Bernardi~\cite{ber12}, Chapuy~\cite{chapuy11},
Chapuy, F\'{e}ray and Fusy~\cite{cff}, Goupil and Schaeffer~\cite{ag98},
Goulden and Nica~\cite{GN05}, Morales and Vassilieva~\cite{MV},
Chen and Reidys~\cite{chr1}, Chen~\cite{chen3} and the references therein.
Though extensively studied,
computing the number of bipartite maps is far from being completely solved.

A novel recursion implying a fundamental identity of Frobenius which enumerates factorizations of a permutation
in group
algebra theory has been discovered by the author recently, and has been used to systematically study the
structure of triple Hurwitz numbers for the first time~\cite{chr-hur}.
Based on that framework, we continue developing this work for the enumeration of bipartite maps and maps. 
Certainly, we also put the study of bipartite maps into the more general algebraic problem
of permutation products or factorizations.
Specifically, our starting point is to count pairs of permutations $(\sigma_1, \sigma_2)$ whose
product has $m$ cycles, where $\sigma_1$ is of cycle-type $\alpha$ and $\sigma_2$ is of cycle-type $\gamma$.
Here we do not require the transitivity of the generated group by $\sigma_1, \sigma_2$.
So, in a sense, ``unconnected" bipartite maps are also considered.
Our main contributions are as follows:
\begin{itemize}
\item[(i)] We present some new general characterizations for the desired number of permutation pairs in terms of symmetric functions;
\item[(ii)] When $\sigma_1$ is limited to a full cycle, we obtain a far-reaching explicit formula of the desired number for any $\gamma$.
The celebrated Harer-Zagier formula~\cite{harer-zagier} (see also in later studies~\cite{zag,GN05,ber12,cff,pittel}) succinctly follows from the formula for $\gamma=[2^n]$, while the
remarkable Jackson's formula~\cite{jac88} (see also in later studies~\cite{sch-vass,ber-morales,cff}) easily follows from summing over all $\gamma$ with the same number of parts;
\item[(iii)] When $\sigma_1$ is limited to a full cycle, in addition to the explicit formula mentioned in (ii), we prove a dimension-reduction
recursion (Theorem~\ref{thm:dim-reduction}) which allows us to compute the number for any $\gamma$ (regardless of the number of parts) from that of the case where $\sigma_2$ is also a full cycle.
The corresponding number of the latter is given by the well-known Zagier-Stanley formula~\cite{stan3,zag};
\item[(iv)] When $\sigma_1$ is limited to a full cycle, we obtain a simple sum formula for $\gamma=[1^t,p,n-t-p]$.
To the best of our knowledge, only the cases $t=0,p=0$ (the Zagier-Stanley formula~\cite{stan3,zag,cms,fv,chr1,chen4}) and $t=1,p=0$ (a consequence of a Boccara's result~\cite{boccara})
are known before.
\end{itemize}

The organization of the paper is as follows.
In Section~\ref{sec2}, we introduce relevant notation of the character theory and present a brief review of the novel recursion
about permutation products.
In Section~\ref{sec3}, we present a general characterization for the number
of arbitrary bipartite maps in terms of Schur symmetric functions.

In the last Section~\ref{sec4}, we focus on one-face bipartite maps.
First, we prove an explicit formula for one-face bipartite maps with any (hyperedge type or {equivalently} white vertex degree distribution) $\gamma$.
We also show one-face bipartite map numbers are essentially determined by a polynomial of the parts of $\gamma$.
Next, we obtain a simple sum formula for $\gamma=[1^t,p,n-t-p]$ for the first time.
After that, we prove a dimension-reduction
recursion which allows us to compute the number for any $\gamma$ (regardless of the number of parts) from that of
$\gamma$ having only one part.
Finally, we show how we can immediately recover the celebrated Harer-Zagier formula and Jackson's formula
which have been the respective main objectives of many classical works in the field.
Direct combinatorial proofs for some of the results in this paper are open and highly desired.

\section{Review of the character theory}\label{sec2}

Denote by $\mathfrak{S}_n$ the symmetric group on $[n]:=\{1,2,\ldots,n\}$.
A permutation $\pi$ can be written as a product of its disjoint cycles, i.e., the orbits of the cyclic group generated
by $\pi$ acting on $[n]$.
The length distribution of the cycles of $\pi$ is called the cycle-type of $\pi$ and denoted by $ct(\pi)$.
One usually writes $ct(\pi)$ as an (integer) partition
of $n$. A partition $\lambda$ of $n$, denoted by $\lambda \vdash n$,
is usually represented by a nonincreasing positive integer sequence $\lambda=(\lambda_1, \lambda_2,\ldots, \lambda_k)$
such that $\lambda_1+\lambda_2+\cdots+\lambda_k=n$. The number $k$ is called the length of $\lambda$, denoted by $\ell(\lambda)$.
Then, $ct(\pi)=\lambda$ means the multiset consisting of the lengths of the cycles of $\pi$ is the same as the multiset $\{\lambda_1, \ldots, \lambda_k\}$.
Another representation of $\lambda$ is in the form $\lambda=[1^{m_1}, 2^{m_2},\ldots, n^{m_n}]$
indicating that there are $m_i$ of $i$'s in the partition.
We often discard the entry $i^{m_i}$ in case $m_i=0$ and write the entry $i^1$ simply as $i$.
Obviously, $\sum_{i=1}^n i m_i=n$ and $\ell(\lambda)=\sum_{i=1}^n m_i$.
One more relevant quantity, denoted by $Aut(\lambda)$, is
the number of permutations of the entries of $\lambda$ that
fix $\lambda$. Clearly, $Aut(\lambda)=\prod_{i=1}^n m_i!$.
In this paper,
we will use the two representations of partitions interchangably, whichever is more convenient.

It is well known that a conjugacy class of $\mathfrak{S}_n$ contains
permutations of the same cycle-type.
So, the conjugacy classes can be indexed by partitions of $n$.
Let $\mathcal{C}_{\lambda}$ denote the one indexed by $\lambda$.
If $\pi \in \mathcal{C}_{\lambda}$, then the number of cycles contained in $\pi$ is $\ell(\pi)=\ell(\lambda)$.
Moreover, the number of elements contained in $\mathcal{C}_{\lambda}$ is well known to be 
$$
|\mathcal{C}_{\lambda}|=\frac{n!}{z_{\lambda}}, \quad \mbox{where $z_{\lambda}= \prod_{i=1}^n i^{m_i} m_i!$.}
$$
Recall the generating function of the signless Stirling
numbers of the first kind $\stirlingI{n}{k}$ is given by:
$$
\sum_{k=1}^n (-1)^{n-k}\stirlingI{n}{k} x^k = x(x-1)\cdots (x-n+1).
$$
The number $\stirlingI{n}{k}$ counts permutations on $[n]$ having exactly $k$ cycles.

From the representation theory of the symmetric group $\mathfrak{S}_n$, we know that
the number of irreducible representations is the same as the number of its conjugacy classes.
Consequently, we can index the irreducible representations by the partitions of $n$ as well.
We write the character associated to the irreducible representation indexed by $\lambda$
as $\chi^{\lambda}$ and the dimension of the irreducible representation as $f^{\lambda}$. Suppose $C$ is a conjugacy class of $\mathfrak{S}_n$ indexed by $\alpha$ and $\pi \in C$.
We view $\chi^{\lambda}(\pi), \, \chi^{\lambda}(\alpha), \, \chi^{\lambda}(C)$ as the same, and we trust the context to prevent confusion.

 Let 
\begin{align*}
\mathfrak{m}_{\lambda,m}= \prod_{u\in \lambda} \frac{m+c(u)}{h(u)}, \qquad 	\mathfrak{c}_{\lambda,m} &= \sum_{d=0}^m (-1)^d {m \choose d} \mathfrak{m}_{\lambda,m-d}.
\end{align*}
where for $u=(i,j)\in \lambda$, i.e., the $(i,j)$ cell in the Young diagram of $\lambda$, $c(u)=j-i$ and $h(u)$ is the hook length of 
the cell $u$. 
The following theorem was proved in~\cite{chr-hur}.

\begin{theorem}[Chen~\cite{chr-hur}] \label{thm:main-recur11} 
	Let $\xi_{n,m}(C_1,\ldots, C_t)$ be the number of tuples $( \sigma_1,\sigma_2,\ldots ,\sigma_t)$ 
	such that the permutation $\pi=  \sigma_1\sigma_2\cdots \sigma_t$ has $m$ cycles, where $\sigma_i $
	belongs to a conjugacy class ${C}_i$. Then we have
	\begin{align}\label{eq:main-recur11}
		\xi_{n,m}(C_1,\ldots, C_t) 
		=\sum_{k = 0}^{n-m} (-1)^k \stirlingI{m+k}{m}   W_{n,m+k}(C_1,\ldots, C_t),
	\end{align}
	where  
	\begin{align}\label{eq:WW}
		W_{n,m+k}(C_1,\ldots, C_t)= \frac{\prod_{i=1}^{t}|C_i|}{(m+k)!} \sum_{\lambda \vdash n} \mathfrak{c}_{\lambda,m+k}	 \big\{f^{\lambda}\big\}^{-t+1} \prod_{i=1}^t \chi^{\lambda}(C_i).
	\end{align}

\end{theorem}

A permutation in $\mathfrak{S}_n$ with only one cycle is called a full cycle or an $n$-cycle,
and a permutation of the cycle-type $[1^j, n-j]$ is called an $(n-j)$-cycle.
The summation in eq.~\eqref{eq:WW} can be simplified if one of the conjugacy classes corresponds to full cycles.
This is presented in the following proposition.

 \begin{proposition}[Chen~\cite{chr-hur}]
	\begin{align}\label{eq:W}
	W_{n,m}(\mathcal{C}_{(n)}, C_1,\ldots, C_t)=\frac{(n-1)!\prod_{i=1}^{t}|C_i|}{m!}\sum_{j=0}^{n-1} (-1)^j \frac{{n-1-j \choose n-m} }{{n-1\choose j}^{t-1}} \prod_{i=1}^t \chi^{[1^j,n-j]}(C_i).
	\end{align}
\end{proposition}

\section{General bipartite maps} \label{sec3}


The enumeration of bipartite maps (not necessarily connected) obviously corresponds to the computation of $\xi_{n,m}(\mathcal{C}_{\alpha}, \mathcal{C}_{\gamma})$.
In the following, we simplify the notation by identifying $\xi_{n,m}(\alpha, \gamma)$ and $\xi_{n,m}(\mathcal{C}_{\alpha}, \mathcal{C}_{\gamma})$.
We first provide a general symmetric function characterization of these numbers.
Let $p_{\lambda}({\bf x})$ and $s_{\lambda}({\bf x})$ denote the power sum symmetric function 
and Schur symmetric function in indeterminates ${\bf x}=x_1, x_2,\ldots$ indexed by $\lambda$, respectively.
We refer to Stanley~\cite{stan-ec2} for precise definitions.

\begin{theorem} Let $n>0$. Then, the following holds:
	\begin{align}\label{eq:schur}
		\sum_{m>0} \sum_{\alpha \vdash n} \sum_{\gamma \vdash n} \frac{1}{(n!)^2} \xi_{n,m}(\alpha, \gamma) p_{\alpha}({\bf x}) p_{\gamma}({\bf y}) z^m=
		\sum_{\lambda \vdash n}  \frac{\mathfrak{m}_{\lambda,z}	} {f^{\lambda}}   s_{\lambda}({\bf x}) s_{\lambda}({\bf y}).
	\end{align}
\end{theorem}

\proof Applying eq.~\eqref{eq:main-recur11}, we have
	\begin{align*}
	&\sum_m \sum_{\alpha} \sum_{\gamma} \frac{1}{(n!)^2} \xi_{n,m}(\alpha, \gamma) p_{\alpha}({\bf x}) p_{\gamma}({\bf y}) z^m\\
	=&\sum_m \sum_{\alpha} \sum_{\gamma} \frac{1}{(n!)^2}
	\sum_{k = 0}^{n-m} (-1)^k \stirlingI{m+k}{m}   \frac{|C_{\alpha}| \cdot |C_{\gamma}|}{(m+k)!} \sum_{\lambda \vdash n} \mathfrak{c}_{\lambda,m+k}	 \big\{f^{\lambda}\big\}^{-1}  \chi^{\lambda}(\alpha) \chi^{\lambda}(\gamma) p_{\alpha}({\bf x}) p_{\gamma}({\bf y}) z^m .
\end{align*}
Since 
\begin{align*}
	\frac{1}{n!}  \sum_{\alpha \vdash n} |C_{\alpha}| \chi^{\lambda}(\alpha) p_{\alpha}({\bf x}) =s_{\lambda}({\bf x}),
\end{align*}
the second last formula equals
\begin{align*}
& \sum_m   
\sum_{k = 0}^{n-m} (-1)^k \stirlingI{m+k}{m}   \frac{1}{(m+k)!} \sum_{\lambda \vdash n} \mathfrak{c}_{\lambda,m+k}	 \big\{f^{\lambda}\big\}^{-1}   s_{\lambda}({\bf x}) s_{\lambda}({\bf y}) z^m \\
=&    
\sum_{j\geq 1} {z \choose j} \sum_{\lambda \vdash n} \mathfrak{c}_{\lambda,j}	 \big\{f^{\lambda}\big\}^{-1}   s_{\lambda}({\bf x}) s_{\lambda}({\bf y}) \\
=&  \sum_{\lambda \vdash n}  \frac{\mathfrak{m}_{\lambda,z}	} {f^{\lambda}}   s_{\lambda}({\bf x}) s_{\lambda}({\bf y}).
\end{align*}
This completes the proof. \qed

Note that $p_{\lambda}({\bf x})$ form a basis of symmetric functions.
In theory, we can obtain $\xi_{n,m}(\alpha, \gamma)$ by first computing the right-hand side of eq.~\eqref{eq:schur}, e.g., using the Jacobi-Trudi determinant formula,
and then taking the coefficient of $p_{\alpha}({\bf x}) p_{\gamma}({\bf y}) z^m$.
This is practical at least for small $n$.

\begin{corollary}
	Let $n>0$. Then, we have
	\begin{align}\label{eq:schur-one}
		\sum_{\alpha \vdash n} \sum_{\gamma \vdash n} \frac{1}{(n!)^2} \xi_{n,1}(\alpha, \gamma) p_{\alpha}({\bf x}) p_{\gamma}({\bf y}) =
\sum_{\lambda \vdash n}   \bigg\{	\sum_{k = 0}^{n-1} (-1)^k \frac{\mathfrak{c}_{\lambda,k+1}}{k+1} \bigg\} \frac{s_{\lambda}({\bf x}) s_{\lambda}({\bf y})}{f^{\lambda}} .
	\end{align}
	\end{corollary}
\proof Analogous to the proof of eq.~\eqref{eq:schur}	. \qed

We remark that a different expression was proved in Morales and Vassilieva~\cite{MV}:
\begin{align}
\label{eq:mor-va}
	\sum_{\alpha \vdash n} \sum_{\gamma \vdash n} \frac{1}{(n!)^2} \xi_{n,1}(\alpha, \gamma) p_{\alpha}({\bf x}) p_{\gamma}({\bf y}) =
	\sum_{\lambda \vdash n}  \sum_{\mu \vdash n} \frac{(n-\ell(\lambda))! (n-\ell(\mu))! }{n!(n+1-\ell(\lambda)-\ell(\mu))!} m_{\lambda}({\bf x})
	m_{\mu}({\bf y}) .
\end{align}
It would be interesting to discuss the connection between the above two expressions.

In the forthcoming section, we will present a systematical study of the cases where $\alpha=[n]$,
resulting in a plethora of new and old results in a unified way.  

\section{One-face bipartite maps and maps} \label{sec4}

A (rooted) one-face bipartite maps (modulo the equivalence) with $n$ edges
is a triple of permutations $(\omega, \sigma, \pi)$ on $[n]$ such that
$\omega=\sigma \pi$ and $\omega$ is a fixed permutation in $\mathcal{C}_{(n)}$, say the cyclic permutation
$(1 \quad 2 \quad \cdots \quad n)$.
If $n=2k$ and $\sigma \in \mathcal{C}_{[2^k]}$, the bipartite map is really a rooted one-face map.
The genus $g$ of the bipartite map is determined by 
$$
2-2g=\ell(\omega )+\ell(\sigma)+ \ell(\pi)-n.
$$

We are now interested in enumerating
all one-face bipartite maps  $(\omega, \sigma, \pi)$ of $n$ edges and genus $g$ where the cycle-type of the
hyperedge $\sigma$ is $\gamma=(\gamma_1,\ldots, \gamma_d)$.
By definition of genus, $\ell(\pi)=m=1-2g+n-d$.
It is more convenient to use $u_{n,m}(\gamma)$ to denote the desired number
and make the genus $g$ implicit in our setting.

As mentioned, this class is believed to be the most simple one.
However, even this class is far from being fully solved, in particular, in terms
of providing explicit computation formulas or algorithms.
We will report our made progress in this regard in the rest of the paper.

The following elementary generating functions will be frequently used throughout this paper.
Moreover, $[x^n]f(x)$ will be used to denote the coefficient of the term $x^n$
in the power series expansion of $f(x)$.

\begin{lemma} For any complex number $\alpha$
	\begin{align}
		\sum_{n\geq 0} {\alpha \choose n} y^n=(1+y)^{\alpha}, \quad \sum_{n\geq 0} {n+\alpha \choose n} y^n =\frac{1}{(1-y)^{\alpha+1}}.
	\end{align}
	\end{lemma}

\subsection{Explicit formulas}
The following proposition is constantly used and can be derived from the Murnaghan-Nakayama rule (see e.g.~Jackson~\cite{jac}).

\begin{proposition}\label{prop:hook-poly}
	Let $\alpha=[1^{a_1},2^{a_2},\ldots,n^{a_n}] \vdash n$. Then,
	\begin{align}
		\sum_{j=0}^{n-1} \chi^{[1^j,n-j]}(\alpha) y^j=(1+y)^{-1} \prod_{i=1}^n \big\{1-(-y)^i\big\}^{a_i}.
	\end{align}
\end{proposition}

Now we are ready to prove a general new explicit formula for one-face bipartite maps.

\begin{theorem}[General explicit formula] \label{thm:hypemap-formula}

Let $\gamma=(\gamma_1,\ldots, \gamma_d) \vdash n$. Then,
	\begin{align}
	\mu_{n,m}(\gamma)
	&= |\mathcal{C}_{\gamma}| \sum_{k=0}^{n-m}   \frac{(-1)^k \stirlingI{m+k}{m}}{(m+k)!} \sum_{k_1+\cdots+k_d=n-m-k+1,\atop k_1>0,\ldots, k_d>0} \prod_{i=1}^d {\gamma_i \choose k_i}. \label{eq:hypermap}
\end{align}
\end{theorem}
\proof 
In the light of Theorem~\ref{thm:main-recur11}, this comes down to
computing the $W$-number for $t=2$, $C_1=\mathcal{C}_{(n)}$, $C_2=\mathcal{C}_{\gamma}$ with $\gamma=(\gamma_1,\gamma_2,\ldots, \gamma_d)$, and we have
\begin{align*}
	W_{n,m+k}&=\frac{(n-1)!|\mathcal{C}_{\gamma}|}{(m+k)!}\sum_{j=0}^{n-1} (-1)^j {n-1-j \choose n-m-k}   \chi^{[1^j, n-j]}(C_2)\\
	&=\frac{(n-1)!|\mathcal{C}_{\gamma}|}{(m+k)!} \sum_{j=0}^{n-1}  [y^{n-m-k}] (1+y)^{n-1-j} (-1)^j \chi^{[1^j, n-j]}(C_2)\\
	&=\frac{(n-1)!|\mathcal{C}_{\gamma}|}{(m+k)!}  [y^{n-m-k+1}] \prod_{i=1}^d \{(1+y)^{\gamma_i} -1\}\\
	&=\frac{(n-1)!|\mathcal{C}_{\gamma}|}{(m+k)!} \sum_{k_1+\cdots+k_d=n-m-k+1,\atop k_1>0,\ldots, k_d>0} \prod_{i=1}^d {\gamma_i \choose k_i}.
\end{align*}
The above third equality follows from Proposition~\ref{prop:hook-poly}.
Note that $\mu_{n,g}(\gamma)$ is $\frac{1}{(n-1)!}$ times
the corresponding $\xi$-number since we fix a long cycle. The rest is clear, and the proof follows. \qed

{ Note that if $\lambda_i=0$ or $n<m$, then the above formula returns zero.}
Immediately from eq.~\eqref{eq:hypermap}, we obtain the following unexpected simple formula
which appeared in Jackson~\cite[Cor.~4.6]{jac}

\begin{corollary}[Genus-zero bipartite maps]
	Let $\gamma=[1^{a_1}, \ldots, n^{a_n}]$. Then, the number of genus zero hypermap with hyper-edge type $\gamma$
	is given by
	\begin{align}
		\frac{n!}{a_1! \cdots a_n! \, [n+1-\ell(\gamma)]!}.
	\end{align}
In particular, the number of genus zero maps of $n$ edges is given by the famous Catalan number $C_n=\frac{1}{n+1} {2n \choose n}$.
\end{corollary}

\proof When $g=0$, the number $m=n+1-d$ where $d= \ell(\gamma)$. As a result, in the summations in eq.~\eqref{eq:hypermap},
only the case $k=0$ and $k_i=1$ matters. Therefore,
$$
\mu_{n,m}(\gamma)=\frac{n!}{1^{a_1} \cdots n^{a_n} a_1!\cdots a_n !} \frac{1}{(n+1-d)!} \gamma_1 \cdots \gamma_d =
\frac{n!}{a_1! \cdots a_n! \, [n+1-d]!}.
$$
 When $\gamma=[2^n]$, the last number reduces to $C_n$,
and the proof follows. \qed

We point out that there is a recursion for $\mu_{n,m}(\gamma)$ in Chen and Reidys~\cite[Cor.~18]{chr1}.
However, that recursion is recursive in both the parameters $g$ (i.e., $m$) and $\gamma$ (through the
refinement relation of partitions).
Moreover, a recursion~\cite[Thm.~17]{chr1} for the numbers refining $\mu_{n,m}(\gamma)$ by tracking the cycle-type of $\pi$ as well (i.e., $\xi_{n,1}(\alpha, \gamma)$) 
is also presented there.
There is also a summation formula for $u_{n,m}(\gamma)$ in Goupil and Schaeffer~\cite[Theorem~$4.1$]{ag98}.
Comparatively speaking, our formula eq.~\eqref{eq:hypermap} may be arguably slightly simpler
and more homogeneous in appearance.
A generating function for the refinement numbers is given in Morales and Vassilieva~\cite{MV}, i.e., eq.~\eqref{eq:mor-va}.

Regarding how $\mu_{n,m} (\gamma)$ depends on the parts $\gamma_i$, we have the following polynomial property.

\begin{theorem}[Polynomiality]
	For fixed $g, \, n$ and $d$, there exists a symmetric polynomial $pol(\gamma_1,\ldots, \gamma_d)$ in $\gamma_1, \gamma_2,\ldots, \gamma_d$ with the highest degree $2g$ such that the weighted number 
$$
\frac{Aut(\gamma)}{n!} \mu_{n,m} (\gamma) = pol(\gamma_1,\ldots, \gamma_d)
$$
 for $m=1-2g+n-d$ and any $\gamma=(\gamma_1,\ldots, \gamma_d) \vdash n$.
	
\end{theorem}

\proof According to eq.~\eqref{eq:hypermap}, we first have
\begin{align*}
	\frac{Aut(\gamma)}{n!} \mu_{n,m}(\gamma)= \frac{1}{\gamma_1 \cdots \gamma_d} \sum_{k=0}^{n-m}   \frac{(-1)^k \stirlingI{m+k}{m}}{(m+k)!} \sum_{k_1+\cdots+k_d=n-m-k+1,\atop k_1>0,\ldots, k_d>0} \prod_{i=1}^d {\gamma_i \choose k_i}.
\end{align*}
Clearly, each summand for $k$ is either zero or a polynomial in $\gamma_1, \ldots, \gamma_d$
where the total degree is $n-m-k+1$ and the power of $\gamma_i$ for any $1\leq i \leq d$ is at least one.
Consequently, $\frac{Aut(\gamma)}{n!} \mu_{n,m} (\gamma)$ is a symmetric polynomial in $\gamma_1, \gamma_2,\ldots, \gamma_d$ with the highest degree
being $(n-m+1)-d=2g$. 
This completes the proof. \qed

Simple formulas for $\mu_{n,m}(\gamma)$ having no complicated summations (e.g., over partitions)
are rare. To the best of our knowledge, only the cases for $\gamma=(n)$ (the Zagier-Stanley formula~\cite{stan3,zag,cms,fv,chr1,chen3}), $\gamma=(n-1,1)$ (a consequence
of a Boccara's result~\cite{boccara}) and $\gamma=[p^n]$ (Jackson~\cite[Thm.~5.4]{jac})
have relatively simple explicit formulas.
However, thanks to our general formula eq.~\eqref{eq:hypermap}, we may derive some new simple explicit formulas. 
We next present two such results as an illustration.

	\begin{theorem}\label{thm:p}
		Let $\gamma=[1^p, n-p]$ for $0\leq p < n$. For $m >0$ and $n-p-m$ being even, the following is true:
		\begin{align}\label{eq:p-power}
			\mu_{n,m}(\gamma)=\frac{2 \stirlingI{n+1-p}{m}}{(n+1-p)!} |\mathcal{C}_{\gamma}|.
		\end{align}
		If $n-p-m$ is odd, $\mu_{n,m}(\gamma)=0$.
	\end{theorem}

\proof  
{
We first compute
\begin{align*}
	 & 1+ \sum_{n } \sum_{m>0}   \sum_{k=0}^{n-m-p}   \frac{(-1)^k \stirlingI{m+k}{m}}{(m+k)!}   {n-p \choose n-m-k+1-p} x^m y^{n+1-p}\\
	&= 1+\sum_{n} \sum_{i\geq 1}   {x \choose i}    {n-p \choose n-i+1-p}  y^{n+1-p}- \sum_{n } {x \choose n+1-p}  y^{n+1-p}\\
	&=1+ \sum_{i\geq 1} {x \choose i} y^i \sum_{n }     {n-p \choose i-1}  y^{n+1-p-i} - (1+y)^{x} \\
	&=1+ \sum_{i\geq 1} {x \choose i} y^i \frac{1}{(1-y)^i} - (1+y)^{x}\\
	&=(1-y)^{-x} - (1+y)^x.
\end{align*}
In the above computation, the sum for $n$ may be over all integers, as long as we assume
${\alpha \choose k}=0$ if $k<0$ or $\alpha < k$.
By construction, the coefficient of the term $x^m y^{n+1-p}$ agrees with $\frac{\mu_{n,m}(\gamma) }{ |\mathcal{C}_{\gamma}|}$
if $n\geq m$ and $n>p \geq 0$. As such, we can proceed as follows. For $n>p$,
}
the coefficient of the term $ y^{n+1-p}$ on the RHS is then
\begin{align*}
{-x \choose n+1-p} (-1)^{n+1-p}-{x \choose n+1-p}.
\end{align*}
The coefficient of the term $x^m$ in the last formula is given by
$$
\frac{\stirlingI{n+1-p}{m}-(-1)^{n+1-p-m} \stirlingI{n+1-p}{m}}{(n+1-p)!}.
$$
Hence, we obtain
{\color{blue}
$$
\mu_{n,m}(\gamma)=\frac{\stirlingI{n+1-p}{m}-(-1)^{n+1-p-m} \stirlingI{n+1-p}{m}}{(n+1-p)!} |\mathcal{C}_{\gamma}| ,
$$
}
completing the proof. \qed

In eq.~\eqref{eq:p-power}, setting $p=0$, we recover the well-known
Zagier-Stanley formula below.

\begin{corollary}[Zagier-Stanley formula~\cite{zag,stan3}] \label{cor:z-s}
	The number of ways of expressing a fixed $n$-cycle as a product of an $n$-cycle and a permutation
	with $m$ cycles is given by 
	$$
	\stirlingI{n+1}{m} \bigg/{n+1 \choose 2}
	$$
	if $n-m$ is even, and $0$ otherwise.
\end{corollary}

Stanley~\cite{stan3} once asked for combinatorial proofs of Corollary~\ref{cor:z-s}.
Various proofs of combinatorial nature have been presented next.
For instance, Cori, Marcus and Schaeffer~\cite{cms}, Chen and Reidys~\cite{chr1}, and Chen~\cite{chen3}.
In particular, a versatile combinatorial approach was presented in Chen~\cite{chen3} to deal
with a variety of problems concerning products of two long cycles.
In F\'{e}ray and Vassilieva~\cite{fv}, a refinement of the problem considered in Corollary~\ref{cor:z-s} was studied, that is,
enumerating the pairs of long cycles whose product has a given cycle-type.
A simple quantitative relation was obtained in~\cite{fv} and was given an alternative simpler proof in~\cite{chen4}.
It would be very interesting to know any combinatorial proofs of our general case eq.~\eqref{eq:p-power}.

{

\begin{theorem}\label{thm:t-p}
	Let $m>0$, and $\gamma=[1^t, p, n-p-t]$ for $p \geq 1$ and $t\geq 0$. Then, $\frac{\mu_{n,m}(\gamma)}{|\mathcal{C}_{\gamma}|}$
	is given by
\begin{align}
 \sum_{j=1}^{n-t} \frac{(-1)^{n-j-t} -(-1)^{j-m} }{j!}  {p \choose n+1-j-t} \stirlingI{j}{m} ,
\end{align}
if $n-m-t$ is odd, and $0$ otherwise.

\end{theorem}

}

\proof 
First, for $n>p+t$, we have
\begin{align*}
\sum_{k_1+\cdots+k_d=n-m-k+1,\atop k_1>0,\ldots, k_d>0} \prod_{i=1}^d {\gamma_i \choose k_i} &={1 \choose 1}^t  \sum_{j=1}^{n-m-k-t} {p \choose j} {n-p-t \choose n-m-k+1-j-t}\\
&=\sum_{j\geq 1} {p \choose j} {n-p-t \choose n-m-k+1-j-t} - {p \choose n-m-k-t+1}.
\end{align*}
{ The above formula does not hold when $n-m-k-t=-1$. It is better to require $n-m-k-t \geq 0$, i.e., $k \leq n-m-t $.

We next compute
\begin{align*}
		&\sum_{n } \sum_{m>0}   \sum_{k=0}^{n-m-t}   \frac{(-1)^k \stirlingI{m+k}{m}}{(m+k)!}   {p \choose n-m-k-t+1} x^m y^{n+1-t}\\
		=& \sum_{n } \sum_{m>0}   \bigg\{ \sum_{i\geq 1}   \frac{(-1)^{i-m} \stirlingI{i}{m}}{i!}   {p \choose n-i-t+1} - \frac{(-1)^{n-t+1-m} \stirlingI{n-t+1}{m}}{(n-t+1)!} \bigg\} x^m y^{n+1-t} \\
		=& \sum_{n }    \sum_{i\geq 1} {x \choose i}     {p \choose n-i-t+1} y^{n+1-t} -  \sum_{n } {x \choose n-t+1} y^{n+1-t}\\
		=&  \sum_{i\geq 1} {x \choose i} y^i \sum_{n } {p \choose n-i-t+1} y^{n+1-t-i}  -  (1+y)^x  \\
		=&  \sum_{i\geq 1} {x \choose i} y^i  (1+y)^p   - (1+y)^x \\
		=&  (1+y)^p \big \{(1+y)^x -1 \big \}  - (1+y)^x.
\end{align*}
Similarly, we have
\begin{align*}
	&  \sum_{n \geq 0} \sum_{m>0}   \sum_{k=0}^{n-m-t}   \frac{(-1)^k \stirlingI{m+k}{m}}{(m+k)!}   \bigg\{ \sum_{j\geq 1} {p \choose j} {n-p-t \choose n-m-k+1-j-t} \bigg\} x^m y^{n+1-t}\\
=	&  \sum_{n \geq 0} \sum_{m>0}   \sum_{i\geq 1}  \frac{(-1)^{i-m} \stirlingI{i}{m}}{i!}   \bigg\{ \sum_{j\geq 1} {p \choose j} {n-p-t \choose n-i+1-j-t} \bigg\} x^m y^{n+1-t}\\
=	&  \sum_{n \geq 0}   \sum_{i\geq 1}  {x \choose i}  y^ i \bigg\{ \sum_{j\geq 1} {p \choose j} y^j {n-p-t \choose n-i+1-j-t} \bigg\}  y^{n+1-t-i-j}\\
=	&    \sum_{i\geq 1}  {x \choose i}  y^ i \sum_{j\geq 1} {p \choose j} y^j \frac{1}{(1-y)^{i+j-p}} \\
=	&    \sum_{i\geq 1}  {x \choose i}  y^ i  \frac{1}{(1-y)^{i}} \big\{ (1-y)^{-p}-1 \big\}  \frac{1}{(1-y)^{-p}} \\
=	&      \big\{ (1-y)^{-x}-1 \big\}  \big\{ 1- {(1-y)^{p}} \big\}   .
\end{align*}
As a result, we obtain
\begin{align*}
	& \sum_{n \geq 0} \sum_{m>0}   \sum_{k=0}^{n-m-t}   \frac{(-1)^k \stirlingI{m+k}{m}}{(m+k)!}   \sum_{j=1}^{n-m-k-t} {p \choose j} {n-p-t \choose n-m-k+1-j-t} x^m y^{n+1-t}\\
=& \big\{ (1-y)^{-x}-1 \big\}  \big\{ 1- {(1-y)^{p}} \big\} -  \big \{(1+y)^x -1 \big \} \big \{(1+y)^p -1 \big \}  +1.
\end{align*}
The coefficient of the term $ y^{n+1-t}$ with $n+1-t>0$ on the RHS of the last formula is 
\begin{align*}
	\sum_{j=1}^{n-t} 	{-x \choose j} (-1)^{j}(-1)^{n-j-t} {p\choose n+1-j-t}-\sum_{j=1}^{n-t} {x \choose j} {p \choose n+1-t-j}.
\end{align*}
The coefficient of the term $x^m$ in the last formula is next obtained as
\begin{align*}
	&\sum_{j=1}^{n-t} {p \choose n+1-j-t} \frac{(-1)^{n-j-t} \stirlingI{j}{m}-(-1)^{j-m} \stirlingI{j}{m}}{j!}  \\
	=& \sum_{j=1}^{n-t} \frac{(-1)^{n-j-t} -(-1)^{j-m} }{j!}  {p \choose n+1-j-t} \stirlingI{j}{m} .
\end{align*}
}
The rest is easy to complete and the proof follows. \qed

\begin{corollary}\label{cor:2parts}
	Let $m>0$ and $\gamma=[p, n-p]$ for $p \geq 1$. Then,
	$\frac{\mu_{n,m}(\gamma) }{|\mathcal{C}_{\gamma}|}=0$ if $n-m$ is even, and
if $n-m$ is odd,
\begin{align}\label{eq:2parts}
	\frac{\mu_{n,m}(\gamma) }{|\mathcal{C}_{\gamma}|} = -2 \sum_{j=m}^n   {p \choose n+1-j} \frac{ s(j,m) }{j!} .
\end{align}
\end{corollary}

In either eq.~\eqref{eq:p-power} or eq.~\eqref{eq:2parts}, setting $p=1$, we obtain
the following formula which can be also obtained by applying a result of Bocarra~\cite{boccara}.
The latter states that there are $2(n-2)!$ ways for writing any odd permutation as a product
of an $n$-cycle and a permutation of $(n-1)$-cycle.

\begin{corollary}
	The number of ways of expressing a fixed $n$-cycle as a product of a permutation in $\mathcal{C}_{(n-1,1)}$ and a permutation
	with $m$ cycles is given by 
	$$
\frac{2}{n-1}	\stirlingI{n}{m} .
	$$
	if $n-m$ is odd, and $0$ otherwise.
\end{corollary}

\subsection{A dimension-reduction recursion}

{
	
	In this part, we prove a dimension-reduction recursion for $\mu_{n,m}(\gamma)$ which
	we believe significant. In the following, we apply the notation:
	\begin{align*}
		\widetilde{\mu}_{n,m}(\gamma) &=\frac{m!}{n! } \mu_{n,m}(\gamma),\\
		\widetilde{S}_{m,i}(l) &= \sum_{j=1}^i {i \choose j} \frac{ (m+j-i)! \, S(l, m+j-i) }{l!} .
	\end{align*}
	Let $\gamma=[1^{a_1}, \ldots, n^{a_n}] \vdash n$ and suppose $a_i>0$.
	We also define
	$$
	\gamma^{\downarrow (i)} : = [1^{a_1}, \ldots, (i-1)^{a_{i-1}}, i^{a_i-1}, (i+1)^{a_{i+1}}, \ldots, n^{a_n}] \vdash n-i.
	$$
	Note that $\gamma^{\downarrow (i)}$ has one less part than $\gamma$.
	Now, we are in a position to present the dimension-reduction formula for one-face bipartite maps.

	\begin{theorem}\label{thm:dim-reduction}
		Let $\gamma=[1^{a_1}, \ldots, n^{a_n}] \vdash n$. Suppose $a_i>0$. Then, 
		\begin{align}
			&\widetilde{\mu}_{n,m}(\gamma)
			=-\sum_{l> m}   \widetilde{S}_{m,1}(l) \widetilde{\mu}_{n,l}(\gamma)  +\frac{1}{i a_i}     \sum_{l> 0}  \widetilde{S}_{m,i}(l) \widetilde{\mu}_{n-i,l }(\gamma^{\downarrow (i)})  .
		\end{align}
	\end{theorem}

	In view of the dimension-reduction formula, all one-face bipartite map numbers are eventually determined by
	the bipartite map numbers of the form $\mu_{n,m}([n])$, i.e., the Zagier-Stanley formula (Cor.~\ref{cor:z-s}).
	This also allows us to compute the one-face bipartite map numbers $\mu_{n,m}(\gamma)$
	for all $n,\, m$ and $\gamma$ in an efficient recursive way:
	construct a database for the cases where $\gamma$ has only one part first, then build
	the database for the cases where $\gamma$ has two parts using the above recursion,
	and for three parts, and so on.
	Walsh and his coauthors have been building such a kind of database, see e.g., Giorgetti and Walsh~\cite{G-walsh},
	and Walsh~\cite{walsh3}.
	Using the dimension-reduction formula for this purpose should be much more efficient than other approaches.

}

\subsection{Jackson's formula}
By summing over all $\gamma$ with the same number of parts, we easily obtain a formula below that is equivalent to the formula of Jackson~\cite{jac88}.
Jackson's formula has been separately proved in Schaeffer and Vassilieva~\cite{sch-vass}, Bernardi and Morales~\cite{ber-morales},
Chapuy, F\'{e}ray and Fusy~\cite{cff}.

\begin{theorem}
	Let $\mu_{n,m}(d) =\sum_{\lambda \vdash n, \ell(\lambda)=d} \mu_{n,m}(\lambda)$.
	Then, we have
\begin{align*}
	 \sum_{m>0} \sum_{d>0} \mu_{n,m}(d) x^m y^d 
	=&\sum_{k \geq 1} {n-1 \choose k-1} {x \choose k}{y+n-k \choose n-k+1}.
\end{align*}
\end{theorem}

\proof First, applying eq.~\eqref{eq:hypermap}, we compute
\begin{align*}
	\sum_{ \gamma \vdash n \atop \ell(\gamma)=d} \mu_{n,m}(\gamma)
	&= \sum_{ \gamma \vdash n \atop \ell(\gamma)=d} |\mathcal{C}_{\gamma}| \sum_{k=0}^{n-m}   \frac{(-1)^k \stirlingI{m+k}{m}}{(m+k)!} \sum_{k_1+\cdots+k_d=n-m-k+1,\atop k_1>0,\ldots, k_d>0} \prod_{i=1}^d {\gamma_i \choose k_i}\\
	&= \sum_{ \gamma \vdash n \atop \ell(\gamma)=d} |\mathcal{C}_{\gamma}| \sum_{k \geq 1}   \frac{(-1)^{k-m} \stirlingI{k}{m}}{k!} \sum_{k_1+\cdots+k_d=n-k+1,\atop k_1>0,\ldots, k_d>0} \prod_{i=1}^d {\gamma_i \choose k_i}\\
	&= \frac{n!}{d!} \sum_{ \gamma_1+\cdots + \gamma_d= n \atop \gamma_i >0} \frac{1}{\gamma_1 \cdots \gamma_d} \sum_{k \geq 1}   \frac{(-1)^{k-m} \stirlingI{k}{m}}{k!} \sum_{k_1+\cdots+k_d=n-k+1,\atop k_1>0,\ldots, k_d>0} \prod_{i=1}^d {\gamma_i \choose k_i}\\
	&= \frac{n!}{d!}   \sum_{k \geq 1}   \frac{(-1)^{k-m} \stirlingI{k}{m}}{k!} \sum_{k_1+\cdots+k_d=n-k+1,\atop k_1>0,\ldots, k_d>0} \frac{1}{k_1 \cdots k_d} \sum_{ \gamma_1+\cdots + \gamma_d= n \atop \gamma_i >0} \prod_{i=1}^d {\gamma_i-1 \choose k_i-1}\\
	&= \frac{n!}{d!}   \sum_{k \geq 1}   \frac{(-1)^{k-m} \stirlingI{k}{m}}{k!} \sum_{k_1+\cdots+k_d=n-k+1,\atop k_1>0,\ldots, k_d>0} \frac{1}{k_1 \cdots k_d} [x^{k-1}] \frac{1}{(1-x)^{k_1+\cdots +k_d}}\\
	&= \frac{n!}{d!}  \sum_{k \geq 1}   \frac{(-1)^{k-m} \stirlingI{k}{m}}{k!}  {n-1 \choose k-1} \frac{d!}{(n-k+1)!} \stirlingI{n-k+1}{d}  .	
\end{align*}
Next, it is straightforward to obtain
\begin{align*}
	& \sum_m \sum_d \mu_{n,m}(d) x^m y^d \\
	=& \sum_m \sum_d    {n!}  \sum_{k \geq 1}   \frac{(-1)^{k-m} \stirlingI{k}{m}}{k!}  {n-1 \choose k-1} \frac{1}{(n-k+1)!} \stirlingI{n-k+1}{d}  x^m y^d  \\
	=&\sum_{k \geq 1} {n-1 \choose k-1} {x \choose k}{y+n-k \choose n-k+1},
\end{align*} 
completing the proof. \qed

\subsection{The Harer-Zagier formula}

In this part, we focus on the case for $\gamma=[2^n]$. From eq.~\eqref{eq:hypermap}, we immediately obtain a new formula
counting the number of one-face maps by genus.

\begin{theorem}[Map formula] 
	For any $n$ and $g$, the number $\xi_{n,g}$ of one-face maps of $n$ edges and genus $g$ is given by
{
	\begin{align}
	\xi_{n,g}
	&=(2n-1)!! \sum_{k\geq 0}  \frac{(-1)^k \stirlingI{m+k}{m}}{(m+k)!} {n \choose m+k-1} 2^{m+k-1} \label{eq:map-formula},
\end{align}
}
	where  $m=n+1-2g$. In particular, $\xi_{n,0}=\frac{1}{n+1} {2n \choose n}$.

\end{theorem}

\proof In this case, the number $m=n+1-2g$, and
\begin{align*}
W_{2n,m+k}=&\frac{(2n-1)!(2n)!}{(m+k)!n!2^n}  [y^{2n-m-k+1}] \{y^2+2y\}^n\\
=& \frac{(2n-1)!(2n)!}{(m+k)!n!2^n} {n\choose x} 2^x
\end{align*}
for $x+2(n-x)=2n-m-k+1$. Solve to obtain $x=k+m-1$,
and the proof follows. \qed

As for the numbers $\xi_{n,g}$, the celebrated Harer-Zagier recurrence~\cite{harer-zagier} provides a three term recurrence
recursive in both $n$ and $g$ while Chapuy~\cite{chapuy11} has obtained the first recurrence recursive only in $g$.
Our recurrence here is apparently also only recursive in $g$.
The first explicit formula for $\xi_{n,g}$ was obtained by Walsh and Lehman~\cite{walsh1}.
Chapuy's recursion can of course yield an explicit formula.
Some relation between these previous formulas, like deriving one from another, and
some new expressions can be found in~\cite{chr2}.
More generally, for $\gamma=[p^n]$, we have
\begin{align*}
	W_{np,m}
	&=\frac{(np-1)!(np)!}{m!n!p^n}  [y^{np-m+1}] \{(1+y)^p-1\}^n\\
	&=\frac{(np-1)!(np)!}{m!n!p^n} \sum_{i=0}^{n} (-1)^i {n\choose i} {p(n-i) \choose np-m+1}.
\end{align*}

In particular, for the extreme case $m=n(p-1)+1$, we have the well-known (generalized) Catalan number expression~\cite{IJ2,jac}.
$$
\mu_{np,m}([p^n])=\frac{1}{(np-1)!} W_{np,m}=\frac{1}{(np-1)!} \frac{(np-1)!(np)!}{(n(p-1)+1)!n!}=\frac{1}{n(p-1)+1} {np \choose n}.
$$

Finally, we present one of the most simple derivations of the celebrated Harer-Zagier formulas
which provide generating functions for the map numbers below.

\begin{corollary}[Harer-Zagier formulas]
	For $n>0$ and $g\geq 0$, we have
	\begin{align}
		\sum_{g \geq 0} \xi_{n,g} x^{n+1-2g} &=(2n-1)!! \sum_{i\geq 1}  {x \choose i} {n \choose i-1} 2^{i-1}, \\
		1+ 2\sum_{n\geq 0} \sum_{g\geq 0}	\frac{\xi_{n,g}}{(2n-1)!!} x^{n+1-2g} y^{n+1} &=\bigg( \frac{1+y}{1-y} \bigg)^x.
	\end{align}
\end{corollary}
\proof Based on eq.~\eqref{eq:map-formula}, we have
\begin{align*}
\sum_{m>0}	\frac{\mu_{n,m}}{(2n-1)!!} x^m &=\sum_{m>0} \sum_{k= 0}^{n-m}  \frac{(-1)^k \stirlingI{m+k}{m}}{(m+k)!} {n \choose m+k-1} 2^{m+k-1} x^m \\
	 &=\sum_{m>0} \sum_{i\geq 1}  \frac{(-1)^{i-m} \stirlingI{i}{m}}{(i)! } {n \choose i-1} 2^{i-1} x^m\\
	 &= \sum_{i\geq 1}  {x \choose i}  {n \choose i-1} 2^{i-1} ,
\end{align*}
whence the first formula.
We proceed to compute as follows.
\begin{align*}
1+ 2\sum_{n\geq 0} \sum_{m>0}	\frac{\mu_{n,m}}{(2n-1)!!} x^m y^{n+1}
	&=1+ \sum_{n \geq 0}  \sum_{i\geq 1} {x \choose i}  {n \choose i-1} 2^{i}  y^{n+1}\\
	&= 1+ \sum_{i\geq 1} {x \choose i} 2^i y^i \sum_{n\geq i-1} {n \choose i-1}  y^{n+1-i}\\
		&= 1+ \sum_{i\geq 1} {x \choose i} 2^i y^i \frac{1}{(1-y)^i}\\
		&= \bigg( \frac{1+y}{1-y} \bigg)^x.
\end{align*}
This completes the proof of the corollary. \qed

{\bf Remark:} This version contains detailed proofs for Theorem 4.5, 4.6 in the published version in Proc. AMS,
and includes the symmetric function description for general bipartite maps and the derivation of Jackson's formula which did not appear in the previous arXiv versions.

\subsection*{Open problems}
Direct combinatorial proofs for
Theorem~\ref{thm:p}, Theorem~\ref{thm:t-p} or Corollary~\ref{cor:2parts}
would be interesting to know.


\section*{Acknowledgements}

The work was partially supported by the Anhui Provincial Natural Science Foundation of China (No.~2208085MA02)
and Overseas Returnee Support Project on Innovation and Entrepreneurship of Anhui Province (No.~11190-46252022001).

{\noindent \bf Declarations of interest:} none.


\end{document}